\documentstyle[11pt]{article}
\begin{document}

\title{Surfaces with flat normal bundle:  an explicit construction}
\author{{\Large Ferapontov E.V.\thanks{
    Present address:
    Fachbereich Mathematik, SFB 288,
    Technische Universit\"at Berlin,
    10623 Berlin,
    Deutschland,\ \ 
    \hbox{e-mail: {\tt fer@sfb288.math.tu-berlin.de}}}}\\
    Centre for Nonlinear Studies\\
    Landau Institute for Theoretical Physics\\
    Academy of Science of Russia, Kosygina 2\\
    117940 Moscow, GSP-1, Russia\\
    e-mail: {\tt fer@landau.ac.ru}}
\date{}
\maketitle

\newtheorem{theorem}{Theorem}

\pagestyle{plain}

\maketitle

\begin{abstract}

An explicit construction of surfaces with flat normal bundle 
in the Euclidean space $E^n$ (unit hypersphere $S^n$) 
in terms of solutions of certain linear system  is proposed. 
In the case of $E^3$ our formulae can be 
viewed as the direct Lie sphere analog of the generalized Weierstrass 
representation of surfaces in conformal geometry or the Lelieuvre representation
of surfaces in the affine 3-space. 

An explicit parametrization of Ribaucour 
congruences of spheres by three solutions of the linear system is obtained.
In view of the classical Lie correspondence between Ribaucour congruences 
and surfaces with flat normal bundle in the Lie quadric in $P^5$
this gives an explicit representation of surfaces with flat normal bundle in
the 4-dimensional space form  of the Lorentzian signature.
Direct projective analog of this construction is the known parametrization of
$W$-congruences by three solutions of the Moutard equation.
Under the Pl\"ucker embedding
$W$-congruences give rise to surfaces with flat normal bundle in the Pl\"ucker
quadric.

Integrable evolutions of surfaces with flat normal bundle and  parallels 
with the theory of nonlocal Hamiltonian operators of hydrodynamic type are 
discussed in the conclusion.

\end{abstract}

\newpage

\section{Introduction}

Let $M^n$ be a submanifold of the Euclidean space $E^{m+n}$ with the radius-vector
$\bf r$ and $m$ pairwise orthogonal unit normals ${\bf n}^\alpha$, the infinitesimal 
displacements of which are governed by the equations
$$
d{\bf n}^\alpha=\omega^\alpha_\beta~{\bf n}^\beta ~~~ mod ~ TM^n.
$$
The 1-forms $\omega^\alpha_\beta$ define connection 
in the normal bundle of
submanifold $M^n$ \cite{Cartan}. 
Submanifolds with flat normal bundle are
characterized by the existence of the normal  frame ${\bf n}^\alpha$ for which
$\omega^\alpha_\beta=0$, so that
$$
d{\bf n}^\alpha\in TM^n ~~~ {\rm for ~ any} ~ \alpha;
$$
such normals are called parallel in the normal bundle. Submanifolds with flat
normal bundle have been extensively investigated in differential geometry:
see e.g. \cite{Akivis}, \cite{Lumiste}, \cite{Terng} and references therein.
Let $u^1,..., u^n$ be a local coordinate system on $M^n$. Introducing the metric
$(d{\bf r}, d{\bf r})=g_{ij}du^idu^j$ and the Weingarten operators 
$(\stackrel{\alpha}{w})^i_j$ by the formulae
$$
\partial_j{\bf n}^\alpha=(\stackrel{\alpha}{w})^i_j~\partial_i\bf r,
$$
$\partial _j=\partial _{u^j},$
one can write down the Gauss-Codazzi equations of a submanifold
with flat normal bundle  in the form
\begin{equation}
\begin{array}{c}
g_{ik}(\stackrel{\alpha}{w})^k_j=g_{jk}(\stackrel{\alpha}{w})^k_i, ~~~~
\nabla _k(\stackrel{\alpha}{w})^i_j=\nabla _j(\stackrel{\alpha}{w})^i_k, \\
\ \\
R^{ij}_{kl}=\sum _{\alpha =1}^m\{(\stackrel{\alpha}{w})^i_k
(\stackrel{\alpha}{w})^j_l-
(\stackrel{\alpha}{w})^j_k(\stackrel{\alpha}{w})^i_l\}
\end{array}
\label{1}
\end{equation}
where $\nabla$ is the covariant differentiation generated by $g_{ij}$ and
$R^{ij}_{kl}=g^{is}R^j_{skl}$ is the curvature tensor. Moreover, the family of
the Weingarten operators is commutative:
\begin{equation}
[\stackrel{\alpha}{w}, \stackrel{\beta}{w}]=0.
\label{2}
\end{equation}
It was observed recently in \cite{Fer3} that the Gauss-Codazzi equations
(\ref{1}), (\ref{2}) of submanifolds with flat normal bundle coincide with
the skew-symmetry conditions and the Jacobi identities for certain nonlocal
Hamiltonian operators of hydrodynamic type (a brief review of this relationship
is included in the Appendix). 
\bigskip

The class of submanifolds with flat normal bundle is invariant under the
following natural transformations:
\bigskip

1. {\bf Conformal transformations} generated by arbitrary translations,
rotations and inversions in $E^{m+n}$. In fact, paper \cite{Akivis} 
treats submanifolds with flat normal bundle as the objects 
of conformal geometry referring them to a special conformal frame of spheres.
\bigskip

2. {\bf Normal shifts} which transform the radius-vector $\bf r$ into
${\bf r}+t^{\alpha}{\bf n}^{\alpha}$; here $t^{\alpha}=const$ and 
${\bf n}^{\alpha}$ are parallel in the normal bundle \cite{Terng}. Submanifolds
related by a normal shift are called parallel.
\bigskip

3. {\bf Gauss map} associating with any submanifold $M^n$ 
with flat normal bundle in $E^{m+n}$ a submanifold 
$\tilde M ^n$ with flat normal bundle
spanned by one of it's parallel normals ${\bf n}^{\alpha}$
in the unit hypersphere $S^{m+n-1}\subset E^{m+n}$.
\bigskip  

4. {\bf Stereographic projection} mapping submanifold $\tilde M ^n\subset
S^{m+n-1}$ into $E^{m+n-1}$.
\bigskip

In the case of hypersurfaces (it should be emphasized 
that any hypersurface of $E^{n+1}$ or
$S^{n+1}$ automatically has flat normal bundle)
transformations 1, 2 generate a finite-dimensional group of contact
transformations known as the Lie sphere group.

Transformations 3, 4 suggest an inductive construction of submanifolds with flat
normal bundle: given a submanifold in $E^{m+n-1}$, we first project it stereographically
into the unit hypersphere $S^{m+n-1}\subset E^{m+n}$ and then reconstruct a submanifold
$M^n$ in $E^{m+n}$ from the given Gauss image (the last step is essentially 
nonunique and requires integration of certain linear system). 
Note that our definition of the Gauss map differs from the standard one,
which associates with a given submanifold the family of it's normal subspaces 
in the Grassmanian.

In the case of surfaces this procedure allows  to construct an arbitrary surface 
with flat normal bundle starting from an orthogonal net on the plane
$E^2$. It requires solving linear equations only.
As will be shown in sect.2, this step-by-step construction can be combined 
into a simple
formula representing an arbitrary surface with flat normal bundle in terms of
solutions of certain linear system.

The commutativity of the family of Weingarten operators implies the
 existence of the net  of curvature lines. The tangents
to the curves of the net are $n$ common eigendirections of 
$\stackrel{\alpha}{w}$.
This net is not necessarily holonomic (that is, a coordinate net): 
isoparametric
submanifolds are the most important examples of submanifolds with flat
normal bundle and  
nonholonomic (in general) net of curvature distributions \cite{Terng3},
\cite{Terng}, \cite{Fer4}.

From now on we consider submanifolds with flat normal bundle and holonomic net
of curvature lines. Submanifolds of this type naturally arise as the
"coordinate" submanifolds of  orthogonal coordinate systems in the 
Euclidean space. Conversely, any submanifold with flat normal bundle and
holonomic net of lines of curvature can be included
(nonuniquely!) in an orthogonal coordinate system. 
Particularly interesting examples of submanifolds with flat normal bundle and 
holonomic net of curvature lines are provided by the embeddings of the
Lobachevski space $L^n$ into $E^{2n+1}$, see e.g. \cite{Aminov},
\cite{Tenenblat}, \cite{Terng2}; here the existence of curvature line 
parametrization was pointed out in \cite{Cartan2}.

In the coordinates
$u^1,..., u^n$ of the lines of curvature the metric $(d{\bf r}, d{\bf r})$ 
and the Weingarten operators $\stackrel{\alpha}{w}$ become diagonal:
$$
(d{\bf r}, d{\bf r})=g_{ii}(du^i)^2, ~~~ (\stackrel{\alpha}{w})^i_j=
(\stackrel{\alpha}{w})^i\delta^i_j
$$
so that the Gauss-Codazzi equations (\ref{1}) assume the form
\begin{equation}
\begin{array}{c}
\partial_j(\stackrel{\alpha}{w})^i=\partial_j \ln \sqrt {g_{ii}}
~ ((\stackrel{\alpha}{w})^j-(\stackrel{\alpha}{w})^i), ~~~ i\ne j \\
\ \\
R^{ij}_{ij}=\sum_{\alpha =1}^m (\stackrel{\alpha}{w})^i(\stackrel{\alpha}{w})^j;
\end{array}
\label{3}
\end{equation}
all other components of the curvature tensor vanish. Introducing the Lame coefficients
$H_i=\sqrt {g_{ii}}$ and the rotation coefficients $\beta_{ij}$ by the
formulae
$$
\partial_iH_j=\beta_{ij}H_i, ~~~ i\ne j,
$$
one can easily check that $(\stackrel{\alpha}{w})^i$ are representable in the form
$$
(\stackrel{\alpha}{w})^i=H^{\alpha}_i/H_i
$$
where $H^{\alpha}_i$ satisfy the same equations as $H_i$:
$$
\partial_iH^{\alpha}_j=\beta_{ij}H^{\alpha}_i.
$$
Since
$$
R^{ij}_{ij}=-\frac{1}{H_iH_j}(\partial_i \beta_{ij}+\partial_j \beta_{ji}
+\sum_{k\ne i, j}\beta_{ki}\beta_{kj}),
$$
equations (\ref{3}) can be rewritten in the form
\begin{equation}
\begin{array}{c}
\partial_iH^{\alpha}_j=\beta_{ij}H^{\alpha}_i \\
\ \\
\partial_k\beta_{ij}=\beta_{ik}\beta_{kj}, ~~~ k\ne i, j \\
\ \\
\partial_i \beta_{ij}+\partial_j \beta_{ji} +
\sum_{k\ne i, j}\beta_{ki}\beta_{kj}+
\sum_{\alpha=1}^mH^{\alpha}_iH^{\alpha}_j=0, ~~~ i\ne j;
\end{array}
\label{4}
\end{equation}
here equations $(\ref{4})_2$ are the compatibility conditions of  $(\ref{4})_1$.
In the case of surfaces they simplify to
\begin{equation}
\begin{array}{c}
\partial_1H^{\alpha}_2=\beta_{12}H^{\alpha}_1, ~~~~
\partial_2H^{\alpha}_1=\beta_{21}H^{\alpha}_2, \\
\ \\
\partial_1 \beta_{12}+\partial_2 \beta_{21} +
\sum_{\alpha=1}^mH^{\alpha}_1H^{\alpha}_2=0;
\end{array}
\label{5}
\end{equation}
in this form they appear in \cite{Eisenhart}, p. 162.

The formulae for submanifolds with flat normal bundle assume much more symmetric form
in the unit hypersphere (the case of the Euclidean space follows by a 
stereographic projection). The advantage of the unit hypersphere is that
unlike the Euclidean case the radius-vector and the unit normals can be treated
on equal footing.
Let ${\bf w}^1$ and ${\bf w}^{\alpha}, ~\alpha=2,..., m$ 
be the radius-vector and the parallel normals of 
$M^n\subset S^{m+n-1}\subset E^{m+n}$,
respectively. We arrange them in an $(m+n)\times m$ matrix $W$ 
which satisfies the equation
$$
W^tW=E_m 
$$
in view of the orthonormality of ${\bf w}^1, ~ {\bf w}^{\alpha}$. The requirement
that $u^1,..., u^n$ are coordinates of lines of curvature implies that the matrices
$\partial_i W$ are of rank one for any $i$. 
\bigskip

In sect.2 we propose a direct construction of surfaces $M^2\subset S^{m+1}
\subset E^{m+2}$ with flat normal bundle by constructing an $(m+2)\times m$
matrix $W$ satisfying the above properties in terms of $m$ arbitrary solutions
$(\kappa^1, s^1), ... , (\kappa^m, s^m)$ of the linear system
\begin{equation}
\begin{array}{c}
\kappa_x=\tan \varphi ~ s_x \\
\kappa_y=-\cot \varphi ~ s_y
\end{array}
\label{rs}
\end{equation}
where $\varphi (x, y)$ is an arbitrary function of $x=u^1, y=u^2$. This 
construction resembles the so-called N-fold Ribaucour transformation known from the theory of
integrable systems (see, e.g., \cite{Liu}) and provides a compact representation
of the iterative procedure discussed above. We prove that an arbitrary surface
with flat normal bundle can be obtained (locally) within our construction.
Particular 
examples of surfaces with flat normal bundle were discussed in
\cite{Palmer} (isothermic surfaces with flat normal bundle), \cite{Vranceanu}
(surfaces flat normal bundle and zero Gaussian curvature); various transformations
of surfaces with flat normal bundle were discussed in \cite{Eisenhart}. The approach
of \cite{Eisenhart} is based on the remark that the radius-vector $\bf r$
of a surface $M^2\subset S^{m+1}\subset E^{m+2}$ with flat normal bundle 
satisfies the Laplace equation
$$
{\bf r}_{xy}=a{\bf r}_x+b{\bf r}_y,
$$
supplemented by a quadratic constraint $({\bf r}, {\bf r})=1$; thence the name:
quadratic conjugate nets.

Construction of sect.2 can be viewed as a direct linearization of
system (\ref{5}) transforming it into the "decoupled" form
$$
\begin{array}{c}
\partial_1H^{\alpha}_2=\beta_{12}H^{\alpha}_1, ~~~~
\partial_2H^{\alpha}_1=\beta_{21}H^{\alpha}_2 \\
\ \\
\partial_1 \beta_{12}+\partial_2 \beta_{21} =0.
\end{array}
$$

Under the change of variables
$\psi^1=-s_y/\sin \varphi, ~~ \psi^2=s_x/\cos \varphi$,
system (\ref{rs}) transforms to a more familiar Dirac operator
$$
\psi^1_x=\varphi_y\psi^2, ~~~
\psi^2_y=-\varphi_x\psi^1;
$$
however, system (\ref{rs}) will be more convenient for our purposes.
\bigskip

In sect.3 we discuss the case of surfaces in $E^3$ and derive the explicit formulae
for the main geometric quantities (such as the radii of principle curvature, fundamental
forms and the Lie-invariant density) in terms of  linear system 
(\ref{rs}). The representation of surfaces in $E^3$ in terms of two solutions
of linear system (\ref{rs}) is a direct analog of the generalized Weierstrass
representation of surfaces in conformal geometry in terms 
of the 2-dimensional Dirac operator \cite{Kon2} and the Lelieuvre representation
of affine surfaces in terms of the Moutard equation \cite{Lelieuvre}.

\bigskip

Sect.4 generalizes construction of surfaces with flat normal bundle to the
case of submanifolds of arbitrary dimension carrying holonomic net 
of curvature lines.
\bigskip

In sect.5 we discuss surfaces with flat normal bundle in the 4-dimensional space of constant
curvature $S^{3, 1}$ of the Lorentzian signature. These surfaces 
can be obtained as 
the images of  Ribaucour congruences of spheres under the Lie sphere map
\cite{Eisenhart}.
In this particular case our construction allows to "parametrize" 
Ribaucour congruences (and hence surfaces with flat normal bundle
in $S^{3, 1}$) by three solutions of system (\ref{rs}).

\bigskip

Construction of sect.5 is the  direct Lie sphere analog of the known
construction in projective differential geometry
parametrizing  W-congruences by three solutions of the Moutard equation.
Under the Pl\"ucker embedding
W-congruences correspond to surfaces with flat normal bundle in the
Pl\"ucker quadric $S^{2, 2}$ \cite{Eisenhart}. We recall this 
construction in sect.6.

\bigskip

Considering linear system (\ref{rs}) as the Lax operator of
the $(2+1)$-dimensional integrable modified Veselov-Novikov (mVN) hierarchy
we define integrable evolutions of surfaces with flat normal bundle
in the spirit of \cite{Kon2}. Although these evolutions are not completely well-defined,
they have a number of interesting properties: in particular, the first
local integral of the mVN hierarchy
$$
\int \int \varphi_x\varphi_y ~ dxdy
$$
coincides with the simplest Lie-invariant functional in  Lie sphere
geometry. This construction is sketched in sect.7.

\bigskip

Parallels between the theory of submanifolds with flat normal bundle and nonlocal
Hamiltonian operators of hydrodynamic type are drawn in the Appendix.

\section{Construction of surfaces with flat normal bundle in a hypersphere}

In this section we propose a construction of surfaces with flat normal bundle 
in terms of solutions of a linear system
$$
\begin{array}{c}
\kappa_x=\lambda (x, y) ~ s_x \\
\kappa_y=\mu (x, y) ~ s_y
\end{array}
$$
where $\lambda, \mu$ are constrained by $\lambda \mu=-1$.
For several reasons it will be convenient to represent this system in the form
(\ref{rs}):
$$
\begin{array}{c}
\kappa_x=\tan \varphi ~ s_x \\
\kappa_y=-\cot \varphi ~ s_y.
\end{array}
$$ 
Let $(\kappa^1, s^1), ... , (\kappa^m, s^m)$ be $m$ arbitrary
solutions of (\ref{rs}). First we introduce two matrices: the $2\times m$
matrix
$$
U=\left(
\begin{array}{c}
\kappa^1~~...~~\kappa^m\\
s^1~~...~~s^m
\end{array}
\right)
$$
and the $m\times m$ matrix $V$ with the elements 
$V^{\alpha \beta}$ defined by the formula
\begin{equation}
dV^{\alpha \beta}=\kappa^{\alpha}d\kappa^{\beta}+s^{\alpha}ds^{\beta}.
\label{Vij}
\end{equation}
The right-hand sides in (\ref{Vij}) are closed in view of (\ref{rs}) so that
$V^{\alpha \beta}$ are 
correctly defined up to additive constants. We restrict these
constants by requiring
\begin{equation}
V^{\alpha \alpha}=\frac{(\kappa^{\alpha})^2+(s^{\alpha})^2+1}{2},
\label{Vii}
\end{equation}
\begin{equation}
V^{\alpha \beta}+V^{\beta \alpha}=\kappa^{\alpha}\kappa^{\beta}+
s^{\alpha}s^{\beta}, ~~~~ \alpha \ne \beta;
\label{Vij+Vji}
\end{equation}
both restrictions are compatible with (\ref{Vij}). In the matrix form
conditions (\ref{Vij}), (\ref{Vii}) and (\ref{Vij+Vji}) can be rewritten as follows:
\begin{equation}
dV=U^tdU
\label{dV}
\end{equation}
\begin{equation}
V+V^t=U^tU+E_m
\label{V+Vt}
\end{equation}
where $E_m$ denotes the $m\times m$ identity matrix.
Let us also point out that
\begin{equation}
\begin{array}{c}
U_x=\left(
\begin{array}{c}
\tan \varphi \\
1
\end{array}
\right)\left( s^1_x ~ ... ~ s^m_x \right) \\
\ \\
U_y=\left(
\begin{array}{c}
-\cot \varphi \\
1
\end{array}
\right)\left( s^1_y ~ ... ~ s^m_y \right) 
\end{array}
\label{dU}
\end{equation}
where the right-hand sides are undestood as the products of the $2\times 1$ and
$1\times m$ matrices. Consequently, $U_x$ and $U_y$ are matrices of rank 1.
Let us introduce finally the $(m+2)\times m$ matrix
\begin{equation}
W=\left(
\begin{array}{c}
UV^{-1} \\
\dotfill \\
V^{-1}-E_m
\end{array}
\right)
\label{W}
\end{equation}
which consists of the $2\times m$ upper submatrix $UV^{-1}$
($V$ is assumed to be invertible) and the $m\times m$ lower submatrix
$V^{-1}-E_m$.

\bigskip

{\bf Lemma 1.} {\it Matrix $W$ satisfies the equation $W^tW=E_m$}.

\bigskip

\centerline{Proof:}

As far as
$$
W^t=\left((V^t)^{-1}U^t ~~ \vdots ~~  (V^t)^{-1}-E_m\right)
$$
we have
$$
\begin{array}{c}
W^tW=(V^t)^{-1}U^tUV^{-1} + ((V^t)^{-1}-E_m)(V^{-1}-E_m)= \\
\ \\
(V^t)^{-1}(V+V^t-E_m)V^{-1} + ((V^t)^{-1}-E_m)(V^{-1}-E_m)=E_m
\end{array}
$$
In this calculation we made use of (\ref{V+Vt}). q.e.d.

\bigskip

Thus,  $m$ columns ${\bf w}^1, ..., {\bf w}^m$ of the matrix $W$ are pairvise 
orthogonal unit vectors in $E^{m+2}$. Let us choose, for instance, 
${\bf w}^1$ as the radius-vector of a surface $M^2$ which by a construction lies
in the unit hypersphere $S^{m+1}\subset E^{m+2}$. We are going to show that the remaining
vectors ${\bf w}^{\alpha}, ~ \alpha = 2, ..., m$ are orthogonal to  $M^2$
and can be interpreted as it's unit normals. The proof is based on the following
\bigskip

{\bf Lemma 2.} {\it Matrices $W_x$ and $W_y$ are of rank one.}

\bigskip

\centerline{Proof:}

Let us demonstrate this for the matrix $W_x$:
$$
\begin{array}{c}
W_x=\left(
\begin{array}{c}
U_xV^{-1}-UV^{-1}V_xV^{-1} \\
\dotfill \\
-V^{-1}V_xV^{-1}
\end{array}
\right)=
\left(
\begin{array}{c}
U_xV^{-1}-UV^{-1}U^tU_xV^{-1} \\
\dotfill \\
-V^{-1}U^tU_xV^{-1}
\end{array}
\right)= \\
\ \\
\left(
\begin{array}{c}
E_2-UV^{-1}U^t \\
\dotfill \\
-V^{-1}U^t
\end{array}
\right)U_xV^{-1}.
\end{array}
$$
Since $U_x$ is of rank one, the rank of $W_x$ cannot exeed one as well. q.e.d.

\bigskip

Hence, ${\bf w}^{\alpha}_x$ are proportional to ${\bf w}^1_x$. 
Similarly, ${\bf w}^{\alpha}_y$ 
are proportional to ${\bf w}^1_y$. As far as $TM^2=span~
\{{\bf w}^1_x, {\bf w}^1_y\}$ and 
${\bf w}^{\alpha}$ are orthogonal to ${\bf w}^{\alpha}_x$  
and ${\bf w}^{\alpha}_y$  
(indeed, ${\bf w}^{\alpha}$ are unit vectors), 
the orthogonality of ${\bf w}^{\alpha}$ 
and $TM^2$ follows directly.
In fact we have proved a stronger result: 
since $d{\bf w}^{\alpha}\subset TM^2$, the surface
$M^2$ automatically has flat normal bundle. Moreover, $x$ and $y$ are
coordinates of the lines of curvature. The metric $(d{\bf w}^1, d{\bf w}^1)$
 and the second fundamental forms $(d{\bf w}^1, d{\bf w}^{\alpha})$ 
of the surface $M^2$ are the elements of the matrix
\begin{equation}
dW^tdW=(X^t)^{-1}dU^tdUX^{-1}
\label{dW}
\end{equation}
(we skip this straightforward calculation). Since the $(\alpha , \beta)$-element
of $dU^tdU$ is of the form
$$
d\kappa^{\alpha}d\kappa^{\beta}+ds^{\alpha}ds^{\beta}=
s^{\alpha}_xs^{\beta}_x~\frac{dx^2}{\cos ^2\varphi}+
s^{\alpha}_ys^{\beta}_y~\frac{dy^2}{\sin ^2\varphi},
$$
the matrix $dW^tdW$ contains no mixed terms $dxdy$. This gives another proof of
the fact that in the coordinates $x, y$ all fundamental forms of $M^2$ are diagonal.
We can formulate the main result of this paper:
\bigskip

{\bf Theorem 1.} {\it For any number $m$ of solutions of linear system (\ref{rs})
each column  of the matrix $W$  can be considered as
the radius-vector of a surface $M^2\subset S^{m+1}\subset E^{m+2}$ with flat normal 
bundle parametrized by coordinates $x, y$ of the lines of curvature.
The remaining columns play the role of $m-1$ normals 
which are parallel in the normal bundle. 
The metric and the second fundamental forms of the surface $M^2$ are given by
(\ref{dW}). 

Conversely, an arbitrary surface $M^2$
with flat normal bundle can be obtained (locally) within this construction.}

\bigskip

\centerline{Proof:}

It remains to prove the last statement. 
This proof is constructive and provides an explicit representation
of the surface $M^2$ in terms of solutions of system (\ref{rs}).
Let $p_0$ be a nongeneric point of 
$M^2\subset S^{m+1}\subset E^{m+2}$ (that is, the net of  curvature lines
is a coordinate net in a neighbourhood of $p_0$). Let ${\bf w}^1$
be the radius-vector of $M^2$ and ${\bf w}^{\alpha}, ~ \alpha =2,..., m$,
the set of parallel normals. We arrange the columns ${\bf w}^1,..., {\bf w}^m$
in an $(m+2)\times m$ matrix $W$ which satisfies the equation
\begin{equation}
W^tW=E_m
\label{ort}
\end{equation}
in view of the orthonormality of ${\bf w}^1,..., {\bf w}^m$. Moreover, one can
always choose coordinates in the ambient space in such a way that 
in the point $p_0$ the matrix $W$ would be of the form
$$
W=\left(
\begin{array}{c}
0 \\
\dotfill \\
E_m
\end{array}
\right).
$$
In a neighbourhood of $p_0$ the matrix $W$ can be represented as follows:
$$
W=\left(
\begin{array}{c}
UV^{-1} \\
\dotfill \\
V^{-1}-E_m
\end{array}
\right)
$$
where $U$ and $V$ are the $2\times m$ and $m\times m$ matrices, respectively.
In the point $p_0$ we have
$U=0, ~~~ V^{-1}=2E_m$ so that $V$ is indeed invertible in a 
neighbourhood of $p_0$. Condition (\ref{ort}) implies (\ref{V+Vt}):
$$
V+V^t=U^tU+E_m.
$$
Let us utilise the fact that $W_x$ and $W_y$ are of rank one
(we recall that $x, y$ are coordinates of the lines of curvature). 
Differentiation of $W$ with respect to $x$ results in
$$
W_x=\left(
\begin{array}{c}
U_x-UV^{-1}V_x \\
\dotfill \\
-V^{-1}V_x
\end{array}
\right)V^{-1}
$$
implying that $V_x$ is of rank one and hence can be represented in the form
$$
V_x=\left(
\begin{array}{c}
a_1 \\
\vdots \\
a_m
\end{array}
\right)\left( \xi^1 ~ ... ~ \xi^m \right)
$$
viewed as the product of $m\times 1$ and $1\times m$ matrices. 
For the same reasons, 
$$
U_x=\left(
\begin{array}{c}
\lambda \\
\beta
\end{array}
\right)\left( \xi^1 ~ ... ~ \xi^m \right).
$$
Differentiation of (\ref{V+Vt}) with respect to $x$
gives
$$
(V_x-U^tU_x)+(V_x-U^tU_x)^t=0,
$$
implying that $V_x-U^tU_x$ is skew-symmetric. On the other hand, it is of rank one.
Since the rank of skew-symmetric matrix is necessarily even, we conclude that
$V_x=U^tU_x$. Similarly, $V_y=U^tU_y$, so that
$$
dV=U^tdU
$$
which coincides with (\ref{dV}). By normalizing $\xi$ we can always represent
$U_x$ in the form
$$
U_x=\left(
\begin{array}{c}
\lambda \\
1
\end{array}
\right)\left( \xi^1 ~ ... ~ \xi^m \right).
$$
Similarly,
$$
U_y=\left(
\begin{array}{c}
\mu \\
1
\end{array}
\right)\left( \eta^1 ~ ... ~ \eta^m \right).
$$
The compatibility conditions $U_{xy}=U_{yx}$ imply $\xi^i=s^i_x, ~ \eta^i=s^i_y$  
so that
$$
U=\left(
\begin{array}{c}
\kappa^1~~...~~\kappa^m\\
s^1~~...~~s^m
\end{array}
\right)
$$
where $\kappa^i_x=\lambda s^i_x, ~ \kappa^i_y=\mu s^i_y$. The condition 
$\lambda \mu =-1$ follows from the compatibility of the equations 
(\ref{dV}) for $V$. q.e.d.

\section{Particular case: surfaces in $E^3$}

Here we give a direct construction of surfaces $M^2\subset E^3$ in terms of two
solutions of linear system (\ref{rs}). Under stereographic projection
$S^3\to E^3$ this construction reduces to that of the preceeding section.
We treat this simplest case separately in order to bring together 
the necessary formulae. 

Let $(\kappa^1, s^1)$ and $(\kappa^2, s^2)$ be two solutions of system (\ref{rs}).
Defining the functions $A$ and $B$ by the formulae
$$
dA=\kappa^1d\kappa^2+s^1ds^2, ~~~~
B=\frac{(\kappa^1)^2+(s^1)^2+1}{2}
$$
we introduce a surface $M^2\subset E^3$ with the radius-vector ${\bf r}$ and the 
unit normal ${\bf n}$:
\begin{equation}
{\bf r}=\left(
\begin{array}{c}
\kappa^2-\kappa^1A/B \\
\ \\
s^2-s^1A/B \\
\ \\
-A/B
\end{array}
\right), ~~~~
{\bf n}=\left(
\begin{array}{c}
\kappa^1/B \\
\ \\
s^1/B \\
\ \\
1/B-1
\end{array}
\right)
\label{r,n}
\end{equation}
(in a somewhat different context these formulae were proposed in
\cite{Fer1}, \cite{Fer2}).
A direct calculation results in the Weingarten equations
\begin{equation}
\begin{array}{c}
{\bf r}_x=\rho^1{\bf n}_x \\
{\bf r}_y=\rho^2{\bf n}_y
\end{array}
\label{Weingarten}
\end{equation}
where $\rho^1, \rho^2$ are the radii of principle curvature of the surface $M^2$:
$$
\rho^1=\lambda^1B-A, ~~~ \rho^2=\lambda^2B-A;
$$
here $\lambda^1={s^2_x}/{s^1_x}, ~ \lambda^2={s^2_y}/{s^1_y}$.
Formulae (\ref{Weingarten}) imply that $x, y$ are coordinates of the lines 
of curvature. By a construction the normal ${\bf n}$ and the 
third fundamental form
$$
(d{\bf n}, d{\bf n})=\frac{(d\kappa^1)^2+(ds^1)^2}{B^2}=
\left(\frac{s^1_x}{B}\right)^2\frac{dx^2}{\cos^2 \varphi}+
\left(\frac{s^1_y}{B}\right)^2\frac{dy^2}{\sin^2 \varphi}
$$
depend only on the first solution $(\kappa^1, s^1)$
of system (\ref{rs}). Varying $(\kappa^2, s^2)$, we reconstruct the full class of 
Combescure-equivalent surfaces. All of them have parallel normals and 
parallel tangents to the lines of curvature in the corresponding points $x, y$;
hence they have one and the same spherical image of the lines of curvature. The results
of sect.2 imply that any surface can be represented (locally) by
formulae (\ref{r,n}) in a neighbourhood of the nonumbilic point.

\bigskip

In the discussion of Lie sphere geometry of surfaces in 3-space
Blaschke \cite{Blaschke} introduced the Lie-invariant functional
which assumes the following form in terms of the radii of principle curvature
$\rho^1, \rho^2$ \cite{Fer}:
\begin{equation}
\int \int \frac{\rho^1_x\rho^2_y}{(\rho^1-\rho^2)^2}~dxdy
\label{Lie}
\end{equation}
As far as $A_x=\lambda^1B_x, ~ A_y=\lambda^2B_y$ we have
$$
 \frac{\rho^1_x\rho^2_y}{(\rho^1-\rho^2)^2}= 
\frac{\lambda^1_x\lambda^2_y}{(\lambda^1-\lambda^2)^2}.
$$
Moreover,
$$
\frac{\lambda^1_x\lambda^2_y}{(\lambda^1-\lambda^2)^2}~dx\wedge dy=
\frac{\lambda^1_y\lambda^2_x}{(\lambda^1-\lambda^2)^2}~dx\wedge dy-
d\left(\frac{d\lambda^2}{\lambda^1-\lambda^2}\right)=
\varphi_x\varphi_y~dx\wedge dy-d\left(\frac{d\lambda^2}{\lambda^1-\lambda^2}\right)
$$
(in this calculation we have used the formula 
$s_{xy}=-\varphi_x\kappa_y-\varphi_y\kappa_x$
which follows from (\ref{rs})). Thus for compact surfaces without umbilic points
(for instance, immersed tori) functional (\ref{Lie}) coincides with
$$
\int\int \varphi_x\varphi_y~dxdy.
$$
As we will see in sect.7, this functional can be interpreted as the first 
conservation law of the (2+1)-dimensional 
integrable mVN hierarchy associated with linear system (\ref{rs}).

\section{Construction of submanifolds with flat normal bundle carrying
coordinate net of curvature lines}

The case of submanifolds $M^n$ of dimension $n$ greater than two requires
certain modifications in the construction of sect.2. 
We start with a Dirac operator
\begin{equation}
\partial _iH_j=\beta_{ij}H_i, ~~~ i, j=1,..., n, ~~~ i\ne j
\label{Lame}
\end{equation}
where the Lame coefficients $H_i$ and the rotation coefficients
$\beta_{ij}$ are functions of $n$ independent variables $u^i, ~ \partial_i=
\partial _{u^i}$. We require that the rotation coefficients satisfy the 
zero curvature conditions
\begin{equation}
\begin{array}{c}
\partial_k\beta_{ij}=\beta_{ik}\beta_{kj}, ~~~ i\ne j\ne k \\
\ \\
\partial_i\beta_{ij}+\partial_j\beta_{ji}+
\sum_{k\ne i, j}^{n}\beta_{ki}\beta_{kj}=0, ~~~ i\ne j.
\end{array}
\label{rotation}
\end{equation}
Equations (\ref{Lame}), (\ref{rotation}) are well-known in the theory of 
$n$-orthogonal curvilinear coordinate systems \cite{Darboux}.
Let us introduce the so-called direction-cosines: $n$ pairwise orthogonal unit vectors
${\bf X}_i=(X_{1i},..., X_{ni})$ satisfying the equations 
\begin{equation}
\begin{array}{c}
\partial_j{\bf X}_i=\beta_{ij}{\bf X}_j \\
\ \\
\partial_i{\bf X}_i=-\sum_{k\ne i}^{n}\beta_{ki}{\bf X}_k
\end{array}
\label{cosine}
\end{equation}
which are compatible in view of (\ref{rotation})
and define the $n\times n$ orthogonal matrix $X_{ji}$.
In order to construct a submanifold $M^n$ with flat normal bundle in the unit
hypersphere $S^{m+n-1}\subset E^{m+n}$ we choose $m$ arbitrary solutions
${\bf H}^{\alpha}=(H^{\alpha}_1,..., H^{\alpha}_n), ~ \alpha=1,..., m$
of the linear system (\ref{Lame}). Defining $m$ vector-functions
${\bf s}^{\alpha}=(s^{\alpha}_1,..., s^{\alpha}_n)$ by the formulae
\begin{equation}
ds^{\alpha}_i=\sum_{k=1}^n X_{ik}H^{\alpha}_kdu^k
\label{flat}
\end{equation}
which are compatible in view of (\ref{rotation}) and (\ref{cosine})
and taking into account the orthogonality of $X_{ji}$, one can 
easily check the identity 
$$
\sum_{i=1}^n(ds^{\alpha}_i)^2=
\sum_{i=1}^n(H^{\alpha}_i)^2(du^i)^2
$$
implying that $s^{\alpha}_i$ are flat coordinates of the flat diagonal
metric 
$$
\sum_{i=1}^n(H^{\alpha}_i)^2(du^i)^2.
$$
Let us introduce two matrices:  the $n\times m$ matrix
$$
U=\left(
\begin{array}{c}
s^1_1~~...~~s^m_1\\
\dotfill \\
s^1_n~~...~~s^m_n
\end{array}
\right)
$$
and the $m\times m$ matrix $V$ with the elements 
$V^{\alpha \beta}$ defined by the formula
\begin{equation}
dV^{\alpha \beta}=({\bf s}^{\alpha}, d{\bf s}^{\beta})=
\sum_{i=1}^ns^{\alpha}_ids^{\beta}_i
\label{Vij2}
\end{equation}
The right-hand sides in (\ref{Vij2}) are closed in view of (\ref{flat}) so that
$V^{\alpha \beta}$ are 
correctly defined up to additive constants. We restrict these
constants by requiring
\begin{equation}
V^{\alpha \alpha}=\frac{({\bf s}^{\alpha}, {\bf s}^{\alpha})+1}{2}=
\frac{\sum_{i=1}^n(s^{\alpha}_i)^2+1}{2},
\label{Vii2}
\end{equation}
\begin{equation}
V^{\alpha \beta}+V^{\beta \alpha}=({\bf s}^{\alpha}, {\bf s}^{\beta})=
\sum_{i=1}^ns^{\alpha}_is^{\beta}_i, ~~~~ \alpha \ne \beta;
\label{Vij+Vji2}
\end{equation}
both restrictions are compatible with (\ref{Vij2}). In the matrix form
conditions (\ref{Vij2})--(\ref{Vij+Vji2}) 
can be rewritten as follows:
$$
dV=U^tdU
$$
$$
V+V^t=U^tU+E_m.
$$
We also point out that as a consequence of (\ref{flat}) all matrices
$\partial_iU$ are of rank one. Let us introduce finally the
$(m+n)\times m$ matrix $W$
$$
W=\left(
\begin{array}{c}
UV^{-1} \\
\dotfill \\
V^{-1}-E_m
\end{array}
\right)
$$
which consists of the $n\times m$ upper submatrix $UV^{-1}$
($V$ is assumed  invertible) and the $m\times m$ lower submatrix
$V^{-1}-E_m$. It is straightforward to check that 
$$
W^tW=E_m.
$$
Thus,  $m$ columns ${\bf w}^1, ..., {\bf w}^m$ of the matrix $W$ are pairwise 
orthogonal unit vectors in $E^{m+n}$. Choosing any one of them 
(say, ${\bf w}^1$) as the radius-vector of a submanifold $M^n\subset
S^{m+n-1}\subset E^{m+n}$ one can show that the remaining
vectors ${\bf w}^{\alpha}, ~ \alpha = 2, ..., m$ are orthogonal to $M^n$
and parallel in the normal bundle. Moreover, $u^i$ are coordinates 
of the lines of curvature. The proofs repeat those from sect.2. 
It is important to emphasise that for $n > 2$ this construction requires the
solution of the nonlinear $n$-wave system (\ref{rotation}) which does not depend
on the codimension of the embedding.

In the case of surfaces equations (\ref{Lame}) are of the form
$$
\partial _1H_2=-\partial_2\varphi ~ H_1, ~~~
\partial _2H_1=\partial_1 \varphi  ~ H_2
$$
where $\beta_{12}=-\partial_2\varphi, ~~ \beta_{21}=\partial_1\varphi$
as a consequence of (\ref{rotation}). It follows from (\ref{cosine})
that
$$
{\bf X}_1=(\cos \varphi, ~ -\sin \varphi), ~~~
{\bf X}_2=(\sin \varphi, ~ \cos \varphi), 
$$
so that equations (\ref{flat}) assume the form
$$
\begin{array}{c}
ds_1=\cos \varphi ~ H_1du^1+ \sin \varphi ~ H_2du^2 \\
\ \\
ds_2=-\sin \varphi ~ H_1du^1+ \cos \varphi ~ H_2du^2 
\end{array}
$$
implying for $s_1, s_2$ the linear system
$$
\begin{array}{c}
\partial _2s_1=\tan \varphi ~ \partial_2 s_2 \\
\ \\
\partial _1s_1=-\cot \varphi ~ \partial_1 s_2 
\end{array}
$$
which coincides with (\ref{rs}). Thus in the case $n=2$ our construction reduces 
to that of sect.2.

\section{Ribaucour congruences and surfaces with flat normal bundle in $S^{3,1}$}

In this section we describe the construction of surfaces with flat normal
bundle in a 4-dimensional space of constant curvature $S^{3, 1}$
of the Lorentzian signature $(+ + + -)$. This construction is based on the notion of 
Ribaucour congruence of spheres which we briefly recall here for a convenience 
of the reader. 

Let $M^2\subset E^3$ be a surface with the radius-vector ${\bf r}$ and the unit
normal ${\bf n}$ satisfying the Weingarten equations
$$
\begin{array}{c}
{\bf r}_x=\rho^1{\bf n}_x \\
{\bf r}_y=\rho^2{\bf n}_y
\end{array}
$$
where $\rho^1, \rho^2$ are the radii of principle curvature.
A sphere $S^2(R)$ of radius $R$ and center ${\bf r}-R{\bf n}$ 
touches $M^2$ at the point ${\bf r}$.
Specifying $R$ as a function of $x, y$ we obtain an arbitrary congruence 
(that is, two-parameter family) of 
spheres tangent to our surface. Besides 
the surface $M^2$ itself, this congruence has another 
envelope $\tilde M^2$.
\bigskip

{\bf Definition.} {\it A congruence of spheres is called the Ribaucour congruence 
if the lines of curvature on $M^2$ and ${\tilde M}^2$ correspond to each other.}
\bigskip

It is known (see e.g. \cite{Eisenhart}, p. 194-196) that the function $R(x, y)$
giving rise to a congruence of Ribaucour is expressible in the form
\begin{equation}
R=\frac{P}{Q}
\label{R}
\end{equation}
where $P$ and $Q$ are solutions of the linear system
\begin{equation}
\begin{array}{c}
P_x=\rho^1Q_x \\
P_y=\rho^2Q_y.
\end{array}
\label{PQ}
\end{equation}
Thus the construction of Ribaucour congruences with the given envelope $M^2$
reduces to the solution of linear system (\ref{PQ}). One can show 
by a direct calculation that the function $R$ defined by (\ref{R}), (\ref{PQ})
satisfies the nonlinear equation
\begin{equation}
R_{xy}=\left(a+{\ln (\rho^1-R)}_y\right)R_x+
\left(b+{\ln (\rho^2-R)}_x\right)R_y
\label{dR}
\end{equation}
where 
$$
a=\frac{\rho^1_y}{\rho^2-\rho^1}, ~~~~ b=\frac{\rho^2_x}{\rho^1-\rho^2}.
$$
Substitution (\ref{R}) linearises the nonlinear equation (\ref{dR}).

\bigskip

With any sphere of radius $R$ and center 
${\bf \xi}=(\xi^1, \xi^2, \xi^3)$ we associate a point 
$$
{\bf Z}=(z^1 : z^2 : z^3 : z^4 : z^5 : z^6)
$$
in projective space $P^5$
with six homogeneous coordinates 
\begin{equation}
\begin{array}{c}
z^1=\xi^1, ~~~ z^2=\xi^2, ~~~ z^3=\xi^3, \\
\ \\ 
z^4=\frac{1-({\bf \xi}, {\bf \xi} )+R^2}{2}, ~~~ 
z^5=\frac{1+({\bf \xi}, {\bf \xi} )-R^2}{2}, ~~~ z^6=R
\end{array}
\label{Z}
\end{equation}
which satisfy the quadratic relation
\begin{equation}
(z^1)^2+(z^2)^2+(z^3)^2+(z^4)^2-(z^5)^2-(z^6)^2=0.
\label{ZZ}
\end{equation}
This construction is Lie's famous correspondence between spheres in $E^3$
and points on the Lie quadric (\ref{ZZ}) in $P^5$, see e.g. \cite{Lie}, 
\cite{Blaschke}, \cite{Cecil}; coordinates $z^i$ are known as hexaspherical coordinates.
 Introducing inhomogeneous coordinates $\tilde z^i=z^i/z^6, ~
i=1,...,5$, we  rewrite (\ref{ZZ}) in the form
$$
(\tilde z^1)^2+(\tilde z^2)^2+(\tilde z^3)^2+(\tilde z^4)^2-(\tilde z^5)^2=1
$$
which represents the 4-dimensional quadric $S^{3, 1}$
of constant curvature $1$ and the signature $(+ + + -)$ in the Lorentzian
space with the metric $(d{\tilde z}^1)^2+(d{\tilde z}^2)^2+(d{\tilde z}^3)^2+
({\tilde z}^4)^2-(d{\tilde z}^5)^2$.
With any congruence of spheres we thus associate a surface in 
$S^{3, 1}$.
\bigskip

{\bf Proposition.} (\cite{Eisenhart}, p. 253)
{\it Ribaucour congruences of spheres correspond to surfaces 
with flat normal bundle in $S^{3, 1}$.}
\bigskip

Indeed, with ${\bf \xi}={\bf r}-R{\bf n}$ and $R=P/Q$ the 6-vector 
${\bf Z}$ defined as in (\ref{Z})
satisfies the Laplace equation
$$
{\bf Z}_{xy}=\left(a+{\ln (\rho^1-R)}_y\right){\bf Z}_x+
\left(b+{\ln (\rho^2-R)}_x\right){\bf Z}_y.
$$
Hence
$$
\tilde {\bf Z}_{xy}=\left(a+\left(\ln \frac{\rho^1-R}{z^6}\right)_y\right)
\tilde {\bf Z}_x+
\left(b+\left(\ln \frac{\rho^2-R}{z^6}\right)_x\right)\tilde {\bf Z}_y
$$
so that all components of the 5-vector $\tilde {\bf Z}$ lying on the quadric
$S^{3, 1}$ satisfy one and the same Laplace equation. This implies that
$\tilde {\bf Z}$ is the radius-vector of a surface with flat normal bundle in
$S^{3, 1}$.

Now we are ready to give a general formula for surfaces with flat normal bundle
in $S^{3, 1}$: for that purpose we construct Ribaucour congruence of spheres in
$E^3$ choosing three solutions $(\kappa^1, s^1), ~ (\kappa^2, s^2), ~ 
(\kappa^3, s^3)$ of linear system (\ref{rs}). First we introduce the functions
$A, B, V^{31}, V^{32}$ as follows:
$$
\begin{array}{c}
dA=\kappa^1d\kappa^2+s^1ds^2, ~~~~ B=\frac{(\kappa^1)^2+(s^1)^2+1}{2}, \\
\ \\
dV^{31}=\kappa^3d\kappa^1+s^3ds^1, ~~~~ dV^{32}=\kappa^3d\kappa^2+s^3ds^2
\end{array}
$$
and define a surface $M^2\subset E^3$ as in sect.3. It's radius-vector ${\bf r}$ 
and the unit normal ${\bf n}$ are the following:
$$
{\bf r}=\left(
\begin{array}{c}
\kappa^2-\kappa^1A/B \\
\ \\
s^2-s^1A/B \\
\ \\
-A/B
\end{array}
\right), ~~~~
{\bf n}=\left(
\begin{array}{c}
\kappa^1/B \\
\ \\
s^1/B \\
\ \\
1/B-1
\end{array}
\right).
$$
The radii of principal curvature $\rho^1, \rho^2$ are given by
$$
\rho^1=\lambda^1B-A, ~~~ \rho^2=\lambda^2B-A
$$
where $\lambda^1={s^2_x}/{s^1_x}, ~ \lambda^2={s^2_y}/{s^1_y}$.
In order to define a Ribaucour congruence of spheres with the given 
envelope $M^2$
we have to specify a pair of functions $P, Q$ satisfying (\ref{PQ}).
As one can verify directly, these can be choosen as follows:
$$
P=V^{32}-V^{31}A/B, ~~~~ Q=V^{31}/B
$$
so that
$$
R=\frac{P}{Q}=\frac{V^{32}}{V^{31}}B-A.
$$
The corresponding congruence of spheres with centers 
${\bf \xi}={\bf r}-R{\bf n}$ and radii $R$ 
generate a surface with flat normal bundle in $S^{3, 1}$ 
with homogeneous hexaspherical coordinates
$$
\begin{array}{c}
z^1=\kappa^2-\kappa^1V^{32}/V^{31}, \\
\ \\
z^2=s^2-s^1V^{32}/V^{31}, \\
\ \\
z^3=(B-1)V^{32}/V^{31}-A, \\
\ \\
z^4=1-{\kappa^2}^2-{s^2}^2+2(\kappa^1\kappa^2+s^1s^2-A)V^{32}/V^{31}, \\
\ \\
z^5=1+{\kappa^2}^2+{s^2}^2+2(A-\kappa^1\kappa^2-s^1s^2)V^{32}/V^{31}, \\
\ \\
z^6=BV^{32}/V^{31}-A.
\end{array}
$$
As in the case of surfaces with flat normal bundle in a hypersphere, 
this construction 
is entirely expressed in terms of solutions of linear system (\ref{rs}).

\section{$W$-congruences and surfaces with flat normal bundle in $S^{2,2}$}

Construction of surfaces with flat normal bundle in 
a 4-dimensional space of constant curvature $S^{2, 2}$ of the signature
$(+ + - -)$ is based on the notion of $W$-congruences of lines which are the direct
projective analogs of Ribaucour congruences of spheres. Under the
Pl\"ucker embedding $W$-congruences correspond to surfaces with flat normal
bundle in the Pl\"ucker quadric in $P^5$. We recall this construction
following \cite{Finikov}, p.139-142.

Let us consider four arbitrary solutions $\xi^1, \xi^2, \xi^3, \xi^4$
of the Moutard equation
\begin{equation}
\xi _{xy}=Q(x, y)\xi
\label{Q}
\end{equation}
and introduce the functions $S^{ij}$ by the formulae
$$
dS^{ij}=(\xi^i_x\xi^j-\xi^j_x\xi^i)dx+(\xi^j_y\xi^i-\xi^i_y\xi^j)dy, ~~~
i\ne j, ~~~ i, j=1,..., 4,
$$
the right-hand sides of which are closed one-forms in view of (\ref{Q}).
Let $M^2$ be a surface in 3-space with the radius-vector
\begin{equation}
{\bf r}=(S^{23}, S^{31}, S^{12}).
\label{Lelieuvre}
\end{equation}
Formula (\ref{Lelieuvre}) is known as the Lelieuvre representation of surfaces 
in 3-space in the asymptotic parametrization $x, y$ \cite{Lelieuvre}.
By a construction
$$
{\bf r}_x={\bf \xi}_x\times {\bf \xi}, ~~~ {\bf r}_y={\bf \xi}\times {\bf \xi}_y
$$
where ${\bf \xi}= (\xi^1, \xi^2, \xi^3)$. In order to construct a $W$-congruence
with the given focal surface $M^2$ we introduce another surface
$\tilde M^2$ with the radius-vector
$$
\tilde {\bf r}={\bf r}+{\tilde {\bf \xi}}\times {\bf \xi}
$$
where
$$
{\tilde {\bf \xi}}=\frac{1}{\xi^4}(S^{41}, S^{42}, S^{43});
$$
vector $\tilde {\bf \xi}$ can be interpreted as 
the Moutard transformation of ${\bf \xi}$ with the help of the fourth
solution $\xi^4$. As pointed out  in \cite{Finikov},  the lines passing through
 ${\bf r}(x, y)$ and ${\tilde {\bf r}}(x, y)$ are tangent
to both surfaces $M^2$ and $\tilde M^2$. Since $x, y$ are asymptotic
coordinates on both surfaces, this congruence is a $W$-congruence with the two
focal surfaces $M^2$ and $\tilde M^2$.
The  homogeneous coordinates of ${\bf r}$ and $\tilde {\bf r}$ are of the form
\begin{equation}
\begin{array}{c}
{\bf r}=(S^{23}, ~ S^{31}, ~ S^{12}, ~ 1), \\
\ \\
\tilde {\bf r}=(S^{23}\xi^4+{S^{42}\xi^3-S^{43}\xi^2}, ~ 
S^{31}\xi^4+{S^{43}\xi^1-S^{41}\xi^3}, ~
S^{12}\xi^4+S^{41}\xi^2-S^{42}\xi^1, ~ 1). \\
\end{array}
\label{rr}
\end{equation}

\bigskip
 
For any two points ${\bf a}$ and
${\bf b}$ in $P^3$ with  homogeneous coordinates
${\bf a}=(a^1 : a^2 : a^3 : a^4), ~~ {\bf b}=(b^1 : b^2 : b^3 : b^4)$  
the Pl\"ucker coordinates of the line $({\bf a}, {\bf b})$ are the six numbers
$$
(p^{12} : p^{13} : p^{14} : p^{23} : p^{42} : p^{34})
$$
($p^{ij}=a^ib^j-a^jb^i$) viewed as homogeneous coordinates of a point in $P^5$
lying on the Pl\"ucker quadric
\begin{equation}
p^{12}p^{34}+p^{13}p^{42}+p^{14}p^{23}=0.
\label{Plucker}
\end{equation}
 In the coordinates
$$
\begin{array}{c}
z^1=\frac{p^{12}+p^{34}}{2}, ~~~ z^4=\frac{p^{12}-p^{34}}{2}, \\
\ \\ 
z^2=\frac{p^{13}+p^{42}}{2}, ~~~ z^5=\frac{p^{13}-p^{42}}{2}, \\
\ \\
z^3=\frac{p^{14}+p^{23}}{2}, ~~~ z^6=\frac{p^{14}-p^{23}}{2}
\end{array}
$$
equation (\ref{Plucker}) assumes the form
$$
(z^1)^2+(z^2)^2+(z^3)^2-(z^4)^2-(z^5)^2-(z^6)^2=0.
$$
Introducing $\tilde z^i=z^i/z^6, ~
i=1,...,5$ we can rewrite this equation in the form
$$
(\tilde z^1)^2+(\tilde z^2)^2+(\tilde z^3)^2-(\tilde z^4)^2-(\tilde z^5)^2=1
$$
which defines the 4-dimensional quadric $S^{2, 2}$
of constant curvature $1$ and the signature $(+ + - -)$ in the Lorentzian
space with the metric $(d{\tilde z}^1)^2+(d{\tilde z}^2)^2+(d{\tilde z}^3)^2-
({\tilde z}^4)^2-(d{\tilde z}^5)^2$.
With any congruence of lines we thus associate a surface in 
$S^{2, 2}$.
\bigskip

{\bf Proposition.} (\cite{Eisenhart}, p. 254)
{\it $W$-congruences of lines correspond to surfaces 
with flat normal bundle in $S^{2, 2}$.}
\bigskip

Applying Pl\"ucker construction to the congruence $({\bf r}, \tilde {\bf r})$
we obtain six homogeneous Pl\"ukker coordinates
$$
p^{12}=S^{23}(S^{43}\xi^1-S^{41}\xi^3)-S^{31}(S^{42}\xi^3-S^{43}\xi^2)
$$
$$
p^{13}=S^{23}(S^{41}\xi^2-S^{42}\xi^1)-S^{12}(S^{42}\xi^3-S^{43}\xi^2)
$$
$$
p^{14}=S^{43}\xi^2-S^{42}\xi^3
$$
$$
p^{23}=S^{31}(S^{41}\xi^2-S^{42}\xi^1)-S^{12}(S^{43}\xi^1-S^{41}\xi^3)
$$
$$
p^{42}=S^{43}\xi^1-S^{41}\xi^3
$$
$$
p^{34}=S^{42}\xi^1-S^{41}\xi^2
$$
which define a surface with flat normal bundle in $S^{2, 2}$ in terms
of four solutions of the Moutard equation.

\section{Integrable evolutions of surfaces with flat normal bundle}

Integrable evolutions of surfaces governed by $(2+1)$-dimensional integrable
equations have been introduced in \cite{Kon2}. The most 
interesting examples include evolution of surfaces in conformal geometry 
based on the generalized Weierstrass representation \cite{Kon2}, \cite{Pinkall},
\cite{Taiman1},  \cite{Taiman3}, and evolution in projective 
geometry based on the Lelieuvre representation of surfaces in 3-space
\cite{Kon2}, \cite{KonPin}. The main idea is that  
linear systems  used to construct a surface (the two-dimensional
Dirac operator in the case of
Weierstrass representation and the Moutard equation
in the Lelieuvre case) are viewed as the Lax operators of the integrable
$(2+1)$-dimensional hierarchies so that the corresponding t-evolutions
act on the induced surfaces. Here we sketch the construction of the third
integrable evolution based on the linear system (\ref{rs}) which is
relevant to Lie sphere geometry.

Linear system (\ref{rs})
$$
\begin{array}{c}
\kappa_x=\tan \varphi ~ s_x \\
\kappa_y=-\cot \varphi ~ s_y
\end{array}
$$
can be supplemented  with the following $t$-evolution of $\kappa$ and $s$
\begin{equation}
\begin{array}{c}
\kappa_t=\kappa_{xxx}-3\cot \varphi ~ \varphi_x(\kappa_{xx}-\varphi_xs_x)+3p\kappa_x \\
\ \\
s_t=s_{xxx}+3\tan \varphi ~ \varphi_x(s_{xx}+\varphi_x\kappa_x)+3ps_x 
\end{array}
\label{rsx}
\end{equation}
which implies the integrable (2+1)-dimensional  equation for $\varphi$:
\begin{equation}
\begin{array}{c}
\varphi_t=\varphi_{xxx}-\varphi_x^3 +3p\varphi_x \\
\ \\
p_y=(\varphi_x\varphi_y)_x. 
\end{array}
\label{varphix}
\end{equation}
Similarly the $\tau$-evolution
\begin{equation}
\begin{array}{c}
\kappa_{\tau}=\kappa_{yyy}+3\tan \varphi ~ \varphi_y(\kappa_{yy}-\varphi_ys_y)+3q\kappa_y \\
\ \\
s_{\tau}=s_{yyy}-3\cot \varphi ~ \varphi_y(s_{yy}+\varphi_y\kappa_y)+3qs_y 
\end{array}
\label{rsy}
\end{equation}
implies the nonlinear equation
\begin{equation}
\begin{array}{c}
\varphi_{\tau}=\varphi_{yyy}-\varphi_y^3 +3q\varphi_y \\
\ \\
q_x=(\varphi_x\varphi_y)_y. 
\end{array}
\label{varphiy}
\end{equation}
Both these $t$- and $\tau$-evolutions are compatible. Their linear combination
\begin{equation}
\begin{array}{c}
\varphi^{'}=\varphi_{xxx}+\varphi_{yyy}-\varphi_x^3-\varphi_y^3 +3p\varphi_x
+3q\varphi_y \\
\ \\
p_y=(\varphi_x\varphi_y)_x\\
\ \\
q_x=(\varphi_x\varphi_y)_y
\end{array}
\label{mVN}
\end{equation}
is known as the $(2+1)$-dimensional potential mKdV equation, 
or the modified Veselov-Novikov (mVN) equation \cite{Kon1}.
Evolution of surfaces in $E^3$ governed by (\ref{mVN})
was also  discussed  in \cite{Schief}.

Under the change of variables
$$
\psi^1=-\frac{s_y}{\sin \varphi}, ~~~ \psi^2=\frac{s_x}{\cos \varphi}
$$
linear system (\ref{rs}) transforms to a more familiar Dirac operator
$$
\begin{array}{c}
\psi^1_x=\varphi_y\psi^2 \\
\ \\
\psi^2_y=-\varphi_x\psi^1
\end{array}
$$
while evolutions (\ref{rsx}) and (\ref{rsy})
assume the forms
\begin{equation}
\begin{array}{c}
\psi^1_t=\psi^1_{xxx}-3\varphi_{xy} \psi^2_x+3\varphi_yp\psi^2 \\
\ \\
\psi^2_t=\psi^2_{xxx}+3(p\psi^2)_x-3\varphi_{xx}\varphi_x\psi^2 \\
\ \\
p_y=(\varphi_x\varphi_y)_x\\
\end{array}
\label{psix}
\end{equation}
and
\begin{equation}
\begin{array}{c}
\psi^1_{\tau}=\psi^1_{yyy}+3(q\psi^1)_y-3\varphi_{yy}\varphi_y\psi^1 \\
\ \\
\psi^2_{\tau}=\psi^2_{yyy}+3\varphi_{xy} \psi^1_y-3\varphi_xq\psi^1 \\
\ \\
q_x=(\varphi_x\varphi_y)_y, \\
\end{array}
\label{psiy}
\end{equation}
respectively. All these evolutions preserve the integral
\begin{equation}
\int \int \varphi_x\varphi_y~dxdy
\label{integral}
\end{equation}
which is the first conservation law in the mVN hierarchy.
Geometric meaning of this functional
in the context of Lie sphere geometry  was clarified in sect.3.
Evolutions (\ref{rsx}) and (\ref{rsy}) 
induce integrable evolutions of surfaces with flat normal bundle. 
Restricting to the case of surfaces in $E^3$ we see that these evolutions 
preserve the Lie-invariant functional (\ref{Lie}) playing a  role similar to that
of the Willmore functional in conformal geometry and the projective area functional
in projective geometry which are invariants of the
evolutions discussed in \cite{Kon2}, \cite{KonPin}, \cite{Pinkall},
\cite{Taiman1}, \cite{Taiman3}.
Thus evolutions of surfaces 
introduced in this section are essentially  Lie-geometric.

{\bf Remark.} Strictly speaking, evolutions of surfaces introduced above are 
not completely well-defined: they depend on the 
particular parametrization of the surface $M^2$ by coordinates $x, y$ of the
lines of curvature as well as on the nonlocalities $p, q$ entering the
equations (\ref{rsx}) and (\ref{rsy}). The only objects which indeed have
an invariant geometric meaning are the integrals of these evolutions
(corresponding to certain Lie-invariant functionals, (\ref{integral}) being
the simplest of them) and their stationary points. It is probably more correct to
speak about foliations of 3-space by one-parameter families of surfaces
 which include the given surface $M^2$.
This point of view was in fact suggested in \cite{Schief}.
The investigation of Lie-geometric properties of these foliations is beyond
the scope of this paper.

\section{Appendix. Nonlocal Hamiltonian operators and
submanifolds with flat normal bundle}

Let us consider an infinite-dimensional phase space of vector functions
$u=\{u^i(x), i=1,\ldots,n\}$, where the Poisson bracket of two functionals
$F=\int f(u, u_x,\ldots)~dx$ and $G=\int g(u, u_x,\ldots )~dx$ is given by
the formula
\begin{equation}
\{I, J\}=\int ~{{\delta F}\over {\delta u^i}}~A^{ij}~
                {{\delta G}\over {\delta u^j}}~dx;
\label{1.1}
\end{equation}
here $A^{ij}$ is an operator of hydrodynamic type
\begin{equation}
A^{ij}=g^{ij}(u)~d+b^{ij}_k(u)~u^k_x,~~~d={d\over {dx}}. 
\label{1.2}
\end{equation}
The theory of such brackets was developed by Dubrovin and Novikov
in \cite{Dubrovin}. 
A fundamental observation was that this theory is essentially
differential-geometric. Indeed, if we take $det~g^{ij}\ne 0$ (such Poisson
brackets are called nondegenerate) and represent $b^{ij}_k$ in the form
$b^{ij}_k=-g^{is}\Gamma ^j_{sk}$, it is not difficult to show that under
point transformations $\tilde u^i= \tilde u^i(u^1,\ldots,u^n)$ the
coefficients $g^{ij}$ transform as components of a type $(2,0)$ tensor,
while $\Gamma ^j_{sk}$  transform as Christoffel symbols of an affine
connection. The condition for the operator (\ref{1.2}) to be Hamiltonian
(i.e., to define a bracket that is skew-symmetric and satisfies the Jacobi
identity) imposes strict constraints on $g^{ij}$ and $\Gamma ^j_{sk}$.

{\bf Theorem 2.} \cite{Dubrovin}
{\it 1. The bracket defined by (\ref{1.1}) and (\ref{1.2}) is skew-symmetric
if and only if the tensor $g^{ij}$ is symmetric (i.e., defines a
pseudo-Riemannian metric) and the connection $\Gamma ^j_{sk}$ is
compatible with the metric: $\nabla_kg^{ij}=0$.

2. The bracket defined by (\ref{1.1}) and (\ref{1.2}) satisfies the 
Jacobi identity if and only if the connection $\Gamma ^j_{sk}$ is symmetric 
and its curvature tensor vanishes.}

In other words, the metric $g_{ij}~du^idu^j$ is flat (here
$g_{ik}g^{kj}=\delta^j_i$), and $\Gamma ^j_{sk}$ are the coefficients of 
the corresponding
Levi-Civita connection. It follows from this that for Hamiltonian operators
of the form (\ref{1.2}) we have an infinite-dimensional analog of the 
Darboux theorem:
in the flat coordinates $g^{ij}=\epsilon ^i\delta ^{ij}~~(\epsilon ^i=\pm 1),
~~ \Gamma ^j_{sk}=0$, and we obtain a particularly simple expression for
$A^{ij}$ with constant coefficients: $A^{ij}=\epsilon ^i\delta^{ij}d$.

If for Hamiltonian we select the hydrodynamic functional $H=\int h(u)dx$,
where density does not explicitely depend on the derivatives 
$u_x, u_{xx}, \ldots$, we obtain a Hamiltonian system of hydrodynamic type
$$u^i_t=A^{ij}{{\delta H}\over {\delta u^j}}= v^i_j(u)~u^j_x,$$
where the matrix $v^i_j$ is given by the formula $v^i_j=\nabla ^i\nabla _j
h$. One may consult the survey \cite{Tsarev} for the necessary information
concerning differential geometry, integrability and applications of
Hamiltonian systems of hydrodynamic type.
 
The
first nonlocal generalization of Hamiltonian operators of hydrodynamic type
was proposed by Mokhov and the author in \cite{Mokhov}:
\begin{equation}
A^{ij}=g^{ij}d-g^{is}\Gamma ^j_{sk}u^k_x+c~u^i_xd^{-1}u^j_x,~~ c=const.
\label{1.3}
\end{equation}
Formally, $A^{ij}$ can be treated as a linear combination of the local
operator (\ref{1.2}) (henceforth, we assume $det~g^{ij}\ne 0$) and the nonlocal
term $u^i_xd^{-1}u^j_x$.

The conditions required for the operator (\ref{1.3}) to be Hamiltonian depend
on the constant $c$ in a nontrivial way.

{\bf Theorem 3.} \cite{Mokhov}
{\it 1. The bracket defined by (\ref{1.1}) and (\ref{1.3}) is skew-symmetric
if and only if the tensor $g^{ij}$ is symmetric 
and the connection $\Gamma ^j_{sk}$ is
compatible with the metric: $\nabla_kg^{ij}=0$.

2. The bracket defined by (\ref{1.1}) and (\ref{1.3}) satisfies the Jacobi identity
if and only if the metric $g_{ij}du^idu^j$ has
constant curvature $c$, ~ i.e., $R^{ij}_{kl}=c(\delta ^i_k\delta ^j_l-
\delta ^j_k\delta ^i_l)$.}

Further generalizations of nonlocal Hamiltonian operators (\ref{1.3})
lean in the direction of modifying the nonlocal "tail":
\begin{equation}
A^{ij}=g^{ij}d-g^{is}\Gamma ^j_{sk}u^k_x+w^i_ku^k_xd^{-1}w^j_lu^l_x.
\label{1.12}
\end{equation}
The conditions required for the operator (\ref{1.12}) to be Hamiltonian impose
certain restrictions on the metric $g^{ij}(u)$, the connection
$\Gamma ^j_{sk}(u)$ and the operator $w^i_j(u)$:

{\bf Theorem 4.} \cite{Fer3}
{\it 1. The bracket defined by (\ref{1.1}) and (\ref{1.12}) is skew-symmetric
if and only if the tensor $g^{ij}$ is symmetric  
and the connection $\Gamma ^j_{sk}$ is
compatible with the metric: $\nabla_kg^{ij}=0$.

2. The bracket defined by (\ref{1.1}) and (\ref{1.12}) satisfies the Jacobi identity
if and only if the connection $\Gamma ^j_{sk}$ is symmetric, and
the metric $g_{ij}$ (with lower indices) and the operator $w^i_j$ satisfy
the Gauss-Codazzi equations:
$$g_{ik}w^k_j=g_{jk}w^k_i,~~~~ \nabla _kw^i_j=\nabla _jw^i_k,$$
$$R^{ij}_{kl}=w^i_kw^j_l-w^j_kw^i_l,~~~~(R^{ij}_{kl}\equiv g^{is}R^j_{skl}).
$$}

In other words, the classical Gauss-Codazzi equations
of hypersurfaces $M^n$ in a pseudo-Euclidean space $E^{n+1}$ are nothing
but the Jacobi identity for the Poisson bracket (\ref{1.1}), (\ref{1.12})! Here the metric
$g_{ij}$ plays the role of the first quadratic form of $M^n$,  
$w^i_j$, the role of the Weingarten operator (shape-operator). If
$M^n$ is a hyperplane in $E^{n+1}$, it's Weingarten operator vanishes and we
obtain the Hamiltonian operator (\ref{1.2}). If $M^n$ is a unit hypersphere, then
it's Weingarten operator $w^i_j=\delta ^i_j$ yields the 
operator (\ref{1.3}) with $c=1$.

Further generalizations involve "lengthening" the nonlocal
tail of the Hamiltonian operator:
\begin{equation}
A^{ij}=g^{ij}d-g^{is}\Gamma ^j_{sk}u^k_x+
\sum _{\alpha =1}^m(\stackrel{\alpha}{w})^i_ku^k_x
d^{-1}(\stackrel{\alpha}{w})^j_lu^l_x.
\label{1.13}
\end{equation}

{\bf Theorem 5.}~ \cite{Fer3}
{\it 1. The bracket defined by (\ref{1.1}) and (\ref{1.13}) is skew-symmetric
if and only if the tensor $g^{ij}$ is symmetric 
and the connection $\Gamma ^j_{sk}$ is
compatible with the metric: $\nabla_kg^{ij}=0$.

2. The bracket defined by (\ref{1.1}) and (\ref{1.13}) satisfies the Jacobi identity
if and only if the connection $\Gamma ^j_{sk}$ is symmetric, and
the metric $g_{ij}$ (with lower indices), and the set of operators
$\stackrel{\alpha}{w}$ 
satisfy the Gauss-Codazzi equations of submanifolds $M^n\subset E^{n+m}$ 
with flat  normal bundle:
$$
g_{ik}(\stackrel{\alpha}{w})^k_j=g_{jk}(\stackrel{\alpha}{w})^k_i,~~~~
\nabla _k(\stackrel{\alpha}{w})^i_j=\nabla _j(\stackrel{\alpha}{w})^i_k,
$$
$$
R^{ij}_{kl}=\sum _{\alpha =1}^N\{(\stackrel{\alpha}{w})^i_k
(\stackrel{\alpha}{w})^j_l-
(\stackrel{\alpha}{w})^j_k(\stackrel{\alpha}{w})^i_l\},
$$
$$
[\stackrel{\alpha}{w}, \stackrel{\beta}{w}]=0.
$$}

 This remarkable correspondence
between nonlocal Hamiltonian operators  and submanifolds $M^n\subset E^{n+m}$
with flat normal bundle
was clarified in \cite{Fer3}, where it was demonstrated that the 
operator (\ref{1.13})
arises as the result of Dirac reduction of the flat operator $\delta ^{IJ}
{d\over {dx}}$, defined in the ambient space $E^{n+m}$ (here $I, J=1,\ldots,
n+m)$, to a submanifold $M^n$.

Thus, the classification of  Hamiltonian operators of the type
(\ref{1.13}) is equivalent to the classification of submanifolds 
with flat normal bundle.
The relevance of nonlocal Hamiltonian operators to the theory of
integrable systems of hydrodynamic type as well as their further
properties and examples are discussed in \cite{Fer3}.

\section{Acknowledgements}

This research was supported by the RFFI grants 96-01-00166,
96-06-80104, INTAS 96-0770 and the Alexander von Humboldt Foundation.
I would like to thank B.G.~Konopelchenko and 
U.~Pinkall for useful discussions.

\end{document}